# Multi-objective optimization for a composite pressure vessel with unequal polar openings

Lyudmyla Rozova[a,b], Bilal Meemary[a], Salim Chaki[a], Mylene Deléglise-Lagardère[a], Dmytro Vasiukov[a]*

Center for Materials and Processes, IMT Nord Europe, Institut Mines-Télécom, Université de Lille, F-59000 Lille, France.

National Technical University "Kharkiv Polytechnic Institute", 2, Kyrpychova str., Kharkiv 61002, Ukraine

*Corresponding author's e-mail: dmytro.vasiukov@imt-nord-europe.fr

**Abstract**

Multi-objective parametric optimization problem is presented for overwrapped composite pressure vessels under internal pressure for storage and heating water. It is solved using the developed iterative optimization algorithm. Optimal values of design parameters for the vessel are obtained by varying the set of parameters for composite layers, such as the thickness of layers and radii of polar openings, which influence the distribution of fiber angles along the vessel. The suggested optimization methodology is based on the mechanical solution for composite vessels and the satisfaction of the main failure criteria. An innovative approach lies in the possibility of using the developed optimization methodology for designing vessels with non-symmetrical filament winding, which have unequal polar openings on the domes. This became possible due to the development of a special numerical mechanical finite element model of a composite vessel. A specific Python program provides the creation of a model and controls the exchange of data between the modules of the iterative optimization process. The numerical model includes the determination of the distribution of fiber angles on the domes and cylindrical part of the vessel as well as changes in layer thicknesses. The optimization problem solution is provided using a Multi-Island Genetic Algorithm, this type of method showed its efficiency for such applications, by allowing to avoid local solutions. Thus, optimal parameters of a composite vessel were found by minimizing composite mass and thickness and maximizing the strain energy. Test solutions using the developed methodology are presented for three types of composite materials to evaluate their possibility for integration into the vessel design model.

1. **Introduction**

The pressure vessels are critical mechanical structures that are generally used to store gases and liquids above atmospheric pressure, under various pressure and temperature conditions, depending on industrial application [1–3]. Of particular interest are Composite Overwrapped Pressure Vessels (COPV), Type IV, with polymer liners, which are popular due to their weight efficiency and good strength properties [4]. The filament winding, using dry fiber passed through liquid resin or preimpregnated tapes, is the most popular manufacturing process to create COPVs. Winding is applied in helical and hoop wraps for cylindrical vessels. Composite vessels could also be used for storing hot water, thus they must withstand temperature changes, pressure, and chlorine corrosion. It offers an alternative to metallic vessels (stainless or enameled steel) conventionally used in domestic hot water boilers. The goal is a durable design with optimal thermal properties and low energy consumption, good performance, recycling, and the possibility of modifications. At the design stage, the study of mechanical and thermal phenomena and selecting optimal composite structure parameters involves formulating and solving an optimization problem [3,5,6].

Structural optimization is crucial for achieving a balance between manufacturing costs and design characteristics like strength, stability, and weight [7–12]. In recent years, several works have been published focusing on the optimization of COPV to reduce weight and cost while enhancing reliability and safety for a variety of applications [2,3,6]. The composite layers design is one of the main factors that directly influence the working pressure and fatigue performances of COPV. Therefore, most of the researches were focused on: 1) the design of the composite shell; 2) the thickness and number of layers; and 3) fiber orientation, which will have a direct impact on the weight and cost of the vessel. Additionally, the purpose can be an increase in its working pressure and service life extension [13–16]. The solution of optimization problems is often based on mechanical or thermomechanical models of vessels [15,17], using finite element (FE) modeling and analytic approach. To create numerical models of composite overwrapped vessels, some commercial FE software such as Abaqus and Ansys can be used with commercial plug-ins, such as Wound Composite Modeler for Abaqus, or self-scripts [3,6,14,16,18,19]. The numerical models are based on axisymmetric or cyclically symmetrical [1,2,14], three-dimensional models, or numerical models of a separate cylindrical part or domes of the vessel [6].

The formulation of the optimization problem as well as choosing the correct objective function, design variables, constraints, and method of solving are the main components of such a study. For example, in work [16], as the objective function, the authors studied the maximum axial buckling load by variating fiber orientation angles. The solutions on maximum burst pressure, by changing fiber angle variation for COPV are also provided in work [2]. Of particular interest is the variety of methods used to solve the optimization problem. Recently, the application of population-based methods for optimization tasks is increasingly in demand [15,16,20,21]. They are based on the biological phenomenon of natural evolution behavior [6]. Thus, to maximize the strength ratio, based on Tsai-Wu criteria for composite layers, the authors of work [20] find the optimal fiber angles and volume by implementing a hybrid method of Genetic Algorithm (GA) with Particle Swarm Optimization (GA-PSO) coupled with local sensitivity analysis. For minimization of the specific strength and the cost factor and satisfying the failure criteria, the Big Bang-Big Crunch (BB-BC) population algorithm is used in work to optimize lay-up configurations and thicknesses of COPV [22]. Methodologies for weight minimization are performed in works [6,21] using GA and Simulated Annealing (SA) methods. The combinations of GA with other techniques are also presented. Authors of work [18] combine GA with direct search to minimize thickness and maximize the stiffness of composite pipes. The combination of Sequential Quadratic Programming (SQP) with GA was applied in [23] to choose the material properties, such as Young's modulus, Bulk modulus, and Poisson's ratios, maximizing the strength ratios for all laminas, varying the fiber angles. The use of population-based methods in combination with neural networks has also shown its effectiveness in the optimization of composite structures [24–26]. They solve the problem of finding the optimal geometric design and stacking sequence, providing weight reduction.

Pressure vessels often take the form of cylinders, spheres, or cones. The vessel optimization, specifically focusing on the dome shape and its composite layer, was explored in [27,28]. The authors proposed optimal meridian profiles for the dome based on non-geodesic trajectories for fibers, aiming to minimize strain energy. The study of optimal dome contour design using the feasible direction method to maximize the shape factor is discussed in [29]. The methodology of vessel design, based on bend-free methodology for superelipsoid of revolution, searching for optimal stiffness by fiber steering was proposed in works [19,30].

Another important geometrical factor in composite vessel design is the presence of unequal polar openings on each of the domes, which leads to changes in the fiber trajectory and fiber angles. Thus, in [31], the trajectory of the fiber is obtained as the solutions in several design sections, for the geodesic and non-geodesic fiber path, with further approximation for the entire vessel model. Special spline interpolation for designing geodesic and semi-geodesic fiber paths on cylinders and cones is used in the work [32]. Testing of the suggested approach on elbows (pipeline fittings) shows good results. Issues of influence of unequal radius of opening on fiber trajectory are considered also for isotensoid vessels in [33]. Such techniques calculate the fiber trajectory based on a mathematical description of the curvature of the mandrel surface, changing its angles of inclination to the vessel axis, for various vessel geometric shapes, convex and concave [34].

Despite the numerous references presented, there is still a lack of a global methodology for the structural multi-objective optimization of non-symmetric pressure vessels with unequal openings, involving both vessel geometry and composite stacking sequence, under mechanical and thermomechanical loading. In this paper, a novel approach for obtaining optimal solutions for non-symmetrical filament-wound structures is proposed. The proposed structural optimization methodology, as a first stage of research, was realized on the base of a numerical mechanical model of a composite vessel, considering both geodesic and non-geodesic fiber trajectories and unequal polar openings on each dome and for each layer. To demonstrate its efficiency, two classical failure criteria were implemented: one is the principal stresses in the fiber direction and the second is the Tsai-Wu quadratic criterion. The methodology was also applied to three different composite materials (carbon fibers, glass fibers, and flax fiber-based thermoplastics) to provide a comprehensive numerical strategy for minimizing the total mass of the hot water storage tanks.

## 2. Composite pressure vessel design

### 2.1. Filament winding parameters

Two types of layers are used to create a filament winding for pressure vessel: hoop layers, which are 90° with the axis of the vessel (in practice, between 85° and 90°), and helical layers, with angles in a range of 0° to less than 90° (effectively between 5° and 85°) [35,36]. The hoop layers are designed to resist tangential forces, created by internal pressure, and provide radial resistance to the structure. They are located on the cylindrical part of the vessel. Helical layers provide the strength of the reservoir domes and are also involved in the perception of axial loads. Each cross layer consists of two helical layers with opposite balanced helical angles in the form of ±α. The fiber orientation angle of helical winding depends on the opening radius, $r_0$ which is the initial radius of the dome and $r$ is the current radius (Fig. 1) [36,37]:

$$\alpha(r) = \arcsin\left(\frac{r_0}{r}\right) \pm \delta \left(\frac{r-r_0}{R-r_0}\right)^n, \tag{1}$$

where $r$ - is the radial distance between the central axe to a point on the dome; $R$ – is the radial distance on the cylinder; $\delta = \alpha_c - \arcsin\left(\frac{r_0}{r}\right)$ - parameter for geodesic winding controlling; $n$ - friction interpolation parameter. Eq. 1 describes two types of fiber trajectories: geodesic in the case when $\delta = 0$ and non-geodesic $\delta \neq 0$.

A geodesic path ensures that fibers are wound along the shortest distance between two points on a curved surface (Fig. 1). This type of trajectory is the simplest since it does not have friction, which prevents the fiber from sliding along the surface of the mandrel. A non-geodesic trajectory is a deviation from the geodesic path in the presence of friction.

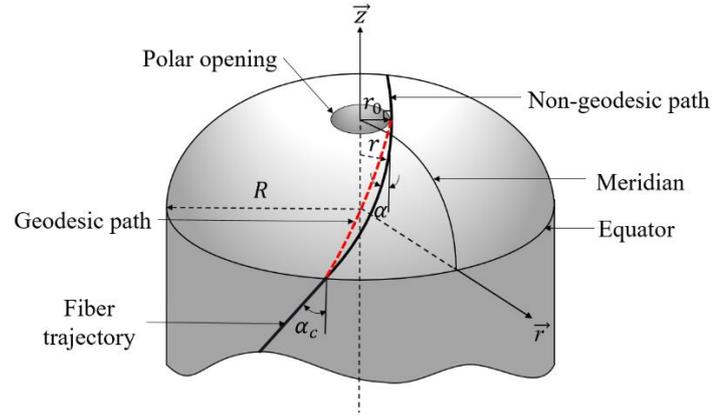

*Fig.1. Fiber path during filament winding.*

The thickness of the composite layer on the vessel depends on the number of rolled-up layers. On the cylindrical part, it is determined by the number of helical and hoop layers. On the domes, it is determined only by the helical layers, that can be wound on this region. Thus, the thickness of the composite changes according to the following law [37]:

$$th(r) = \frac{th_c \cos(\alpha_c)}{\cos(\alpha(r))} \frac{r}{r + 2w_B \left(\frac{R-r}{R-r_0}\right)^4} \; , \qquad (2)$$

where $th_c$ - helical layers thickness on cylinder; $\alpha(r)$ - fiber angle on the dome; $w_B$ - bandwidth. The thickness of the composite layer on the dome also depends on $r_0$. The main originality of this work lies in the consideration that each composite layer has its polar opening. By varying the opening radii, it is possible to ensure the required strength characteristics of the vessel and minimize the mass of a composite structure [37]. On the other hand, the presence of unequal polar openings on different domes of the vessel requires providing a solution that considers changes in the fiber angle on the cylindrical part of the vessel for each composite layer [31,32]. Thus, to complete a global optimization methodology, a special solution should also be provided for the slippage angle of the cylindrical section. In this work, it is assumed that the fiber angle changes linearly along the cylindrical part.

Schematic trajectories of fibers for helical windings of the first and the second layers are shown in Fig. 2a,b respectively. The openings of the first helical layer $r_{0up_1}$ and $r_{0d_1}$ (Fig. 2a) are technological constraints and, thus, the whole surface should be covered by composite. All subsequent helical layers are the layers with openings $r_{0up_i}$ and $r_{0d_i}$ (*i=2...n*, for *n* is the total number of helical layers).

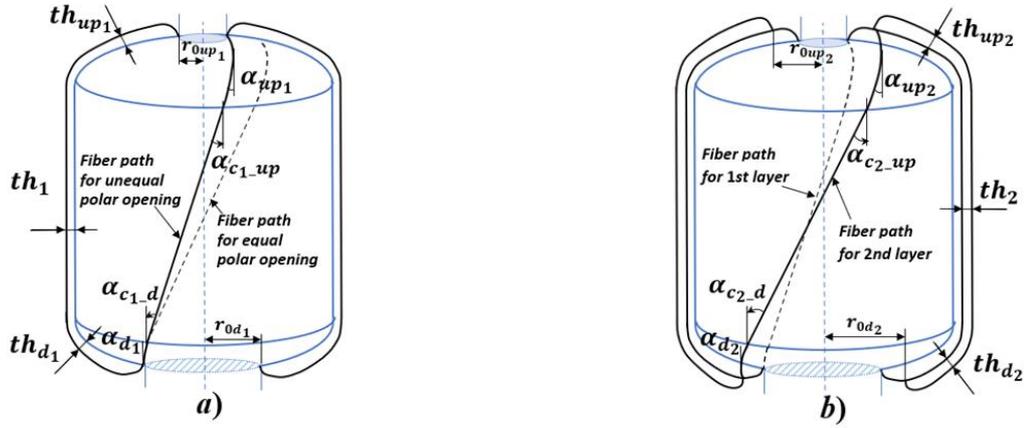

*Fig .2. Schematic representation of helical winding fiber path with unequal openings: a) for one layer; b) for two layers.*

### 2.2. Failure criteria

The composite material considered for COPV can be modeled as the unidirectional composite material with orthotropic or transverse isotropic material behavior [38]. Design criteria can be based on a variety of currently existing failure criteria. In this work, we selected one of the often-used Tsai-Wu quadratic criterion [20]. According to this criterion, the destruction of an orthotropic material is achieved when the following condition is satisfied:

$$F_1\sigma_1 + F_2\sigma_2 + F_{11}\sigma_1^2 + F_{22}\sigma_2^2 + F_{66}\sigma_6^2 + 2F_{12}\sigma_1\sigma_2 = 1 \qquad (3)$$

with, $F_1 = \frac{1}{X_t} - \frac{1}{X_t}$; $F_2 = \frac{1}{Y_t} - \frac{1}{Y_c}$; $F_{11} = \frac{1}{X_tX_c}$; $F_{22} = \frac{1}{Y_tY_c}$; $F_{66} = \frac{1}{S_{12}^2}$; $F_{12} = -\frac{1}{2}\sqrt{F_{11}F_{22}}$,

where $X_t, X_c$ - longitudinal tensile and compressive strength of composite; $Y_t, Y_c$ - transverse tensile and compressive strength of composite; $S_{12}$- shear strength. However, it is classified as an independent criterion that does not consider the effects of the different failure modes (fiber failure; matrix failure; delamination damages) [1,39]. Tsai-Wu criterion can be considered more conservative and may give an overestimated result for the first-ply failure pressure in composite, therefore, for a more accurate determination of failure pressure and failure modes, it is proposed to use dependent criteria, such as Hashin or three-dimensional invariant-based failure criteria [1,40,41]. The purpose of this study was to develop an optimization methodology that allows determining the optimal parameters of the considered structure, where the failure criteria are used as constraints, one of them is the Tsai-Wu criterion, which takes into account the 3-dimensional stress state of the structure. And the other one, as an alternative approach for a quick estimation of the thickness of the plies, is the criterion of maximum principal stresses in the fiber direction, which is used for a "netting analysis", that allows a comparison of the obtained results [37,42]:

$$\sigma_1 < X_t \tag{4}$$

These two criteria were selected for this study, but the methodology is not limited to them.

### 2.3. Material properties

For our study, the following three types of composite materials were used for optimization:

a) Polypropylene (PP), reinforced by glass fibers (GF) (Polystrand, Avient [43])

b) Polyamide (PA), reinforced by carbon fibers (CF) (Toray Cetex [44])

c) Polylactide (PLA) with flax fibers (FF) (Lincore, Depestele [45])

All these materials are available in the market. The collected from the manufacturer and the estimated mechanical properties of those materials are given in Table 1. To complete the missing material elastic constants for transversally isotropic material and strength properties for (Eq. 1), the well-established system of Chamis' rules [38,50] (see Appendix 1, Eq. A1) was used. In this work, we used material data presented in the literature [46–48]. It should be noted that we use the same material polypropylene (PP) for the liner and the composite of the pressure vessel for each solution [46]. For the accuracy of further comparison of calculated results, all values of material properties were recalculated for a 0.45 fiber volume fraction for all composite materials presented in Table 1.

*Table 1. Mechanical properties*

| Material | PP (liner) | GF-PP | CF-PA | FF-PLA |
|---|---|---|---|---|
| $V^f$ | | 0.45 | 0.45 | 0.45 |
| $E_{11}$, $MPa$ | 1300 | 34500 | 90296 | 30643 |
| $E_{22}$, $MPa$ | | 3450 | 6111* | 5621* |
| $\nu_{12}$ | 0.42 | 0.28 | 0.28* | 0.41* |
| $\nu_{13}$, $\nu_{23}$ | | 0.25 | 0.31* | |
| $G_{12}$, $MPa$ | | 2140 | 2763* | 1908* |
| $G_{23}$, $MPa$ | | | 2279* | |
| $X_T$, $MPa$ | 25 | 765 | 1710 | 295 |
| $X_C$, $MPa$ | | 357 | 337 | 170* |
| $Y_T$, $MPa$ | 30 | 13.1 (24*) | 65.8* | 35* |
| $Y_C$, $MPa$ | | 38.3 (36.6*) | 71.6* | 52* |
| $S_{12}$, $MPa$ | | 22.1 | 39 | 40* |
| $\rho$, $10^{-9}$ $t/mm^3$ | | 1.63 | 1.47 | 1.3 |

*values, estimated using Chamis' rules (see Appendix 1 Eq. A1).

### 2.4. Material section assignments

The changes in the fiber angles and the thicknesses of each layer on each dome are considered, according to Eqs. 4 and 5. They depend on the polar openings for each layer and on the upper ($r_{0up_i}$, $\alpha_{up_{ij}}$) and down domes ($r_{0d_i}$, $\alpha_{d_{ij}}$) separately. It should be noted that each composite layer, in the numerical mechanical model, means the set of layers of winding with the same polar opening on the dome, further called a layer. Also, because of the modeling of helical winding, it assumes the existence of paired layers with the angles $\pm\alpha_{up_{ij}}$ and $\pm\alpha_{d_{ij}}$.

On the other hand, for the considered model of a composite vessel, due to the design features, the presence of unequal polar openings during winding on each dome is assumed. As mentioned earlier, the angle of fiber trajectory along the cylindrical part is determined, depending on the value of the polar opening. So, for each dome at the point of intersection of the meridian of the dome with the cylinder, we obtain different fiber angles. This leads to the changes of fiber angles on the cylindrical part of the vessel. The changes in fiber angle in the composite layer are modeled as a piecewise linear approach. To do this, the cylindrical part of the vessel is partitioned into sections (Fig.3*a*). Special arrays of fiber angle values for each composite layer ($\alpha_{c_{ij}}$), depending on the polar radii of the opening on the domes, and for each section are determined.

### 2.5. Geometry, loadings and boundary conditions

In this work, an industrial application example, a new design of a composite water storage tank, was studied. The optimization problem (Fig. 4) was implemented for a composite pressure vessel with a volume of 180 L, the main geometrical parameters and values of applied loads are presented in Table 2.

*Table 2. Main initial values of geometrical parameters for composite vessel and loading conditions*

| Geometrical parameters, mm | | | | | Loads, MPa | |
|---|---|---|---|---|---|---|
| Height of cylinder, $H_c$ | The radius of the cylinder, $R$ | Height of upper dome, $h_{up}$ | Height of down dome, $h_d$ | Liner's thickness, $th_{ln}$ | $P_m$ | $P_{hyd}$ |
| 819 | 250 | 123 | 123 | 4 | 2 | 0.01 |

To find the optimal design of the considered composite vessel, an FE mechanical numerical model was developed within Abaqus/Standard [11] using a specially developed Python script, which covers all geometrical and material features discussed above. Quadrilateral conventional shell elements S4R were used for the cylindrical part and domes (Fig. 3c). The modeling of composite layers was accomplished using the composite section technique for conventional shells in Abaqus, employing the calculated values of angles and thicknesses for different layers. The mesh convergence study determined an average optimal finite element size of 10 mm.

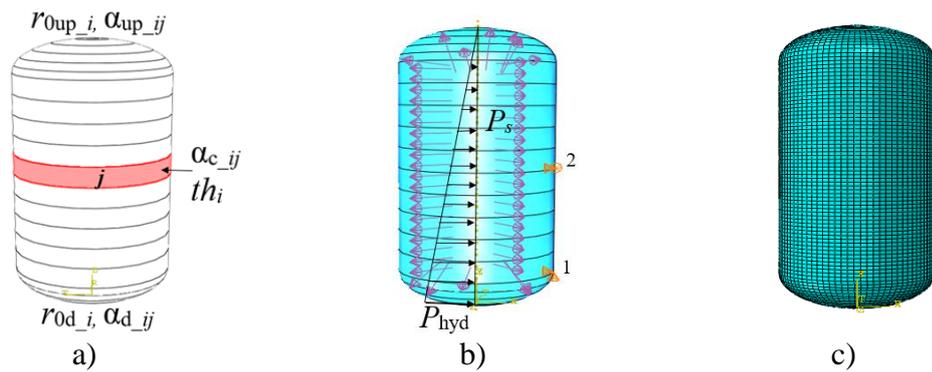

*Fig.3. Numerical model: a) geometrical partitions; b) loadings and boundary conditions; c) finite element model.*

The vessel is subjected to internal uniform pressure (Fig. 3*b*) of amplitude **$P_m$** = 2 MPa. The hydrostatic pressure from water acts on the inner surface of the vessel. This pressure is proportional to gravity and increases linearly with height due to the increasing weight of water, which exerts increasing force on the walls. Therefore, it is zero at the top of the reservoir and maximum at the bottom: **$P_{hyd}$** = $\rho g h$, where $\rho$ - water density; $g$ - acceleration of gravity; $h$ - total height of vessel. The vessel is fixed at two points 1 and 2 in (Fig. 3*b*) were created as the boundary conditions, first, to prevent the vessel from moving in the axial direction (Y: $Uy=0$). In addition to eliminate vessel rotation, the fixations were added in a (X-Z) plane, parallel to the cylinder axis ($Ux=0$ in point 1,2, $Uz=0$ in point 2).

After each simulation, the stress components, failure indices and strain energy were passed to the optimization module.

## 3. Multi-objective optimization problem.

The mathematical formulation of the optimization problem is to find the minimum or maximum of objective function or functions $f_k(\bar{x})$, $k = \overline{(1, l)}$, subjected to constraints $\varphi_j(\bar{x}) \geq 0$, $j =$ $\overline{1,m_1}$; $\psi_j(\bar{x}) = 0$, $j = \overline{m_1 + 1, m}$, so to find the vector of variables $\bar{x}^\star = [x_1^\star, x_2^\star, ..., x_n^\star]^T$ for which these conditions can be satisfied [7,8]:

$$\text{Min}(max)\ f_k(\bar{x}), k = \overline{1,l};$$

$$\underline{\text{Subjected to:}} \tag{5}$$

$$\begin{cases} \varphi_j(\bar{x}) \geq 0, & j = \overline{1, m_1}; \\ \psi_j(\bar{x}) = 0, & j = \overline{m_1 + 1, m}; \end{cases}$$

In this article, we address the challenge through parametric optimization. This approach enables us to identify the optimal combination of parameters that are essential for pressure vessel design. The optimization process is constructed to ensure mechanical performance along with minimization of the total mass, and, thus, directly a total cost of the composite. Two key parameters have a direct influence on the mass: the thickness of each ply and the radius of the polar openings for each helical layer. These openings, in turn, impact the fiber angles along the vessel during winding, subsequently affecting the mechanical response of the vessel.

Thus, the multi-objective parametric optimization problem is formulated as shown in Fig. 4. The objective functions are the total mass of the composite and the total thickness of the composite on the cylindrical part of the vessel. The aim is to minimize these functions to reduce the cost of the studied structure. Simultaneously, we search to maximize strain energy to enhance the mechanical capability of the vessel.

The solution to the optimization task is intended to satisfy two types of failure criteria, serving as constraints. The first criterion is based on the maximum principal stresses in the fiber direction (Eq. 4). This research aims to estimate the thickness for a specific case in the optimization problem and compare the results with values obtained using the estimated formula for the thickness of a composite cylinder in an equal stress state [37,42]. This type of estimation is commonly used for industrial purposes to determine vessel thickness. The second type of solution to the optimization problem is based on the Tsai-Wu failure criterion Eq. 3. Also, as for constraints, manufacturing constraints are used, such as constant polar opening for the first helical layer on each dome.

The summary of the parameters of the optimization problem is listed in Fig. 4. The developed optimization methodology allows us to use as variables the geometrical sizes of vessels, material properties, thicknesses and polar openings for each wounded composite layer.

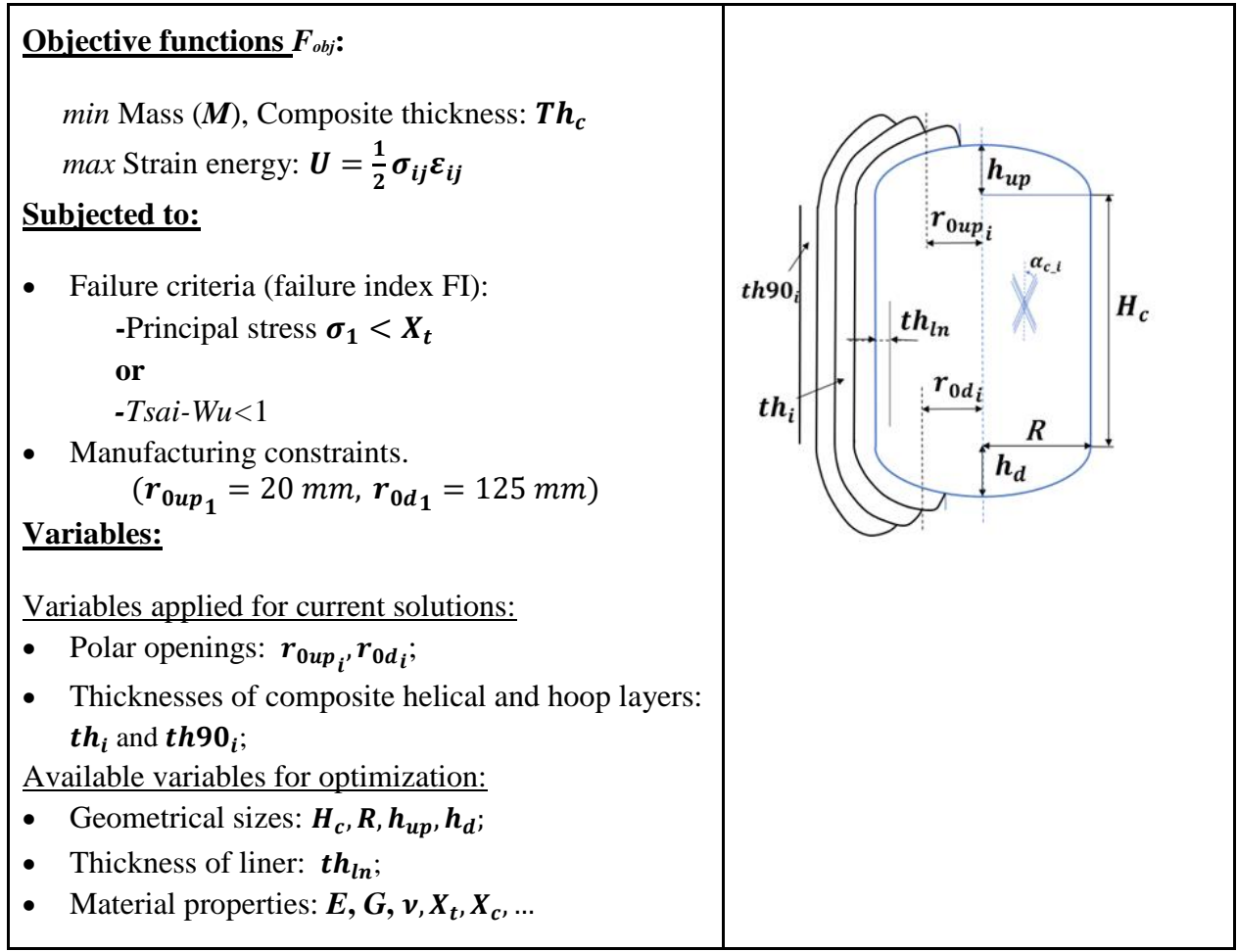

**Objective functions $F_{obj}$:**

  *min* Mass ($M$), Composite thickness: $Th_c$
  *max* Strain energy: $U = \frac{1}{2}\sigma_{ij}\varepsilon_{ij}$

**Subjected to:**

- Failure criteria (failure index FI):
    -Principal stress $\sigma_1 < X_t$
    **or**
    -*Tsai-Wu*<1
- Manufacturing constraints.
    ($r_{0up_1} = 20\ mm$, $r_{0d_1} = 125\ mm$)

**Variables:**

Variables applied for current solutions:
- Polar openings: $r_{0up_i}, r_{0d_i}$;
- Thicknesses of composite helical and hoop layers: $th_i$ and $th90_i$;

Available variables for optimization:
- Geometrical sizes: $H_c, R, h_{up}, h_d$;
- Thickness of liner: $th_{ln}$;
- Material properties: $E, G, \nu, X_t, X_c, ...$

*Fig.4. Summary of the multi-objective parametrical model*

## 4. Optimization method and implemented algorithm

As a numerical platform for the optimization problem, the Simulia Isight [49] software package is used for the automation of all steps of the optimization. To solve mechanical problems Abaqus/Standard FE solver is used. There are several optimization techniques available for parametric optimization in the Isight optimization module. A study of the use of various optimization methods showed that gradient methods are not suitable for the considered problem and have bad convergence. The main cause is that this type of method provides solutions for very local extremums, requires the calculation of the objective functions' gradients at each step, and strongly depends on the initial approximation [8,15,16].

Better results give the application of non-gradient methods from the family of direct search methods, such as Hooke-Jeeves. This technique uses a combination of objective and constraint penalties as the objective functions and does not use derivatives of objective functions. Instead of this, the algorithm examines points near the current point by perturbing design variables, one axis at a time, until an improved point is found. In such a way the favorable direction is found until no more design improvement is possible [7,10]. The disadvantage of such methods: they stop the search when a local minimum or maximum is found, but with a correctly selected initial approximation, a suitable optimal solution can be obtained. The main advantage of this method includes its fast convergence; a small number of iterations are required. The runtime of the solution process for the considered model takes an average of 110 min for 200 iterations (calculations were carried out on a PC with processor Intel(R) Xeon(R) W-1250 CPU @ 3.30GHz and memory RAM 64Gb). Considering the aforementioned, and according to the literature this method can be used for the assessment of the optimal solution [18].

Another type of optimization technique based on GA is also tested for the solution. The provided literature review showed the efficiency of the application of GA for such type of structural parametric multi-objective optimization problems [6,16,18,20,21]. They belong to the family of population methods that simulate the phenomenon of natural evolution. Each design point is an individual for objective function and constraints. An individual with a better value of the objective function and constraints has a higher fitness value. Each population of individuals (a set of design points) is obtained by the genetic operations of "selection", "crossover", and "mutation". A new population of designs is selected from the original set of designs, the best individuals from the previous generations. These techniques are well-suited for different design spaces, and for multi-objective optimization, no need for the calculation of function gradients. However, a long computational time is required [49]. For the presented research, the Multi-Island Genetic

Algorithm was applied (MIGA). Like other genetic algorithms, in MIGA each design point is perceived as an individual. Each individual is represented by a chromosome in which the values of design variables are converted into a binary string of 0 and 1 characters, and participate in forming a population.

Having studied all the stages for designing the considered composite vessel for heating and storing water, a special algorithm for parametric optimization was suggested (Fig. 5).

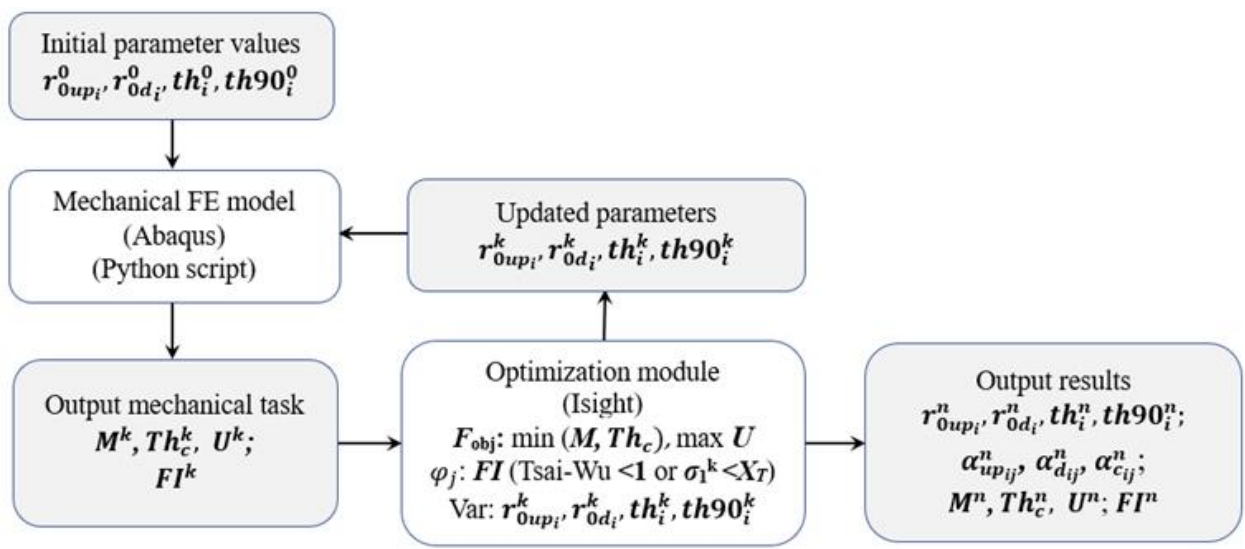

*Fig.5. Proposed optimization algorithm.*

In the proposed algorithm, at each iteration, a finite element model of the vessel is built, using the Python script. Depending on the current values of the polar opening radii $r_{0up_i}, r_{0d_i}$ and the thicknesses of the composite layers ($th_i^k, th90_i^k$) the special arrays of distribution of fiber angles along the vessel ($\alpha_{up_{ij}}, \alpha_{d_{ij}}, \alpha_{c_{ij}}$) are calculated, and also the changes in thicknesses on each dome ($th_{d\,i}^k, th_{up\,i}^k$), according to (Eq. 2). The mass of the composite and its thickness are also determined ($M^k, Th_{c\,i}^k$). The solution to the mechanical problem provides values of strain energy ($U^k$) and failure index ($FI^k$) at each iteration.

It is important to ensure the correct and fast exchange of optimization parameters at each cycle and the calculation of the in-line values of the studied quantities. Data exchange between modules is provided by the Python program using input and output files. Values for objective functions, constraints and parameters are transferred to the optimization module. After the optimization solution, data is exchanged again with the mechanical module. The automatic update of the input and output parameters is also controlled by the Python program.

## 5. Results and discussion

To find the optimal solution according to the proposed algorithm we fix the geometrical sizes of the vessel and vary the polar radii of the openings on each dome ($r_{0up_i}, r_{0d_i}$) for helical composite layers, and the thicknesses of each composite layer ($th_i^k, th90_i^k$). The solution model of the vessel includes four hoop layers and three helical layers, which consist of paired layers with the angles $\pm\alpha_{ij}$. At the same time, the values for the first composite layer remained constant ($r_{0up_1} = 20\ mm$, $r_{0d_1} = 125\ mm$), following the imposed manufacturing features (Fig. 2).

Two groups of solutions are provided: 1) the first group is based on the maximum principal stresses in the fiber direction failure criterion, and 2) the second group is based on the Tsai-Wu quadratic failure criterion. The goal of this work was to develop an efficient methodology for structural parametric optimization of a non-symmetrical composite pressure vessel designed for storing and heating water, aiming to create a functional and intelligent design for the vessel as a whole. The obtained optimal values of the main vessel parameters can guide the further choice of composite material. To confirm the proper performance of the suggested optimization scenario, the first group of solutions, using the maximum principal stresses in the fiber direction failure criterion, serves to estimate and compare the results with the "netting analysis" [37,42]. Such a simplified approach is often used in the classical design of symmetrical composite pressure vessels for quick estimation of the main parameters.

The results of the first group solutions are presented in Table 3 for three types of composite materials presented in Table 1. It gives us a rough estimation of vessel thickness, considering the equal uniform stress state [37,42]. In this case, the total thickness of the composite (both hoop and helical layers) on the cylindrical zone can be evaluated as:

$$th = \frac{3PR}{2\sigma_1} \quad (6)$$

The thickness of the composite in this case directly depends on the value of $X_t$, as the maximum value for $\sigma_1$ which is the stress value in the fiber direction. So, it is the relevant first estimation of vessel thickness, using only strength properties in the fiber direction. The results of the calculations show a good agreement with the values, obtained by Eq. 6, which can prove the developed approach.

*Table 3. The summary of the calculations based on the maximum principal stress criterion*

| Materials | Netting analysis | Max Principal stress | | |
|---|---|---|---|---|
| | Total thickness, according to (Eq. 6), (mm) | Total thickness, obtained (mm) | Difference, (%) | The total mass of composite, (kg) |
| **GF-PP** | 0.98 | 1.04 | 6.1 | 3.3 |
| **CF-PA** | 0.43 | 0.47 | 9.3 | 1.5 |
| **FF-PLA** | 2.54 | 2.65 | 4.3 | 5.9 |

The second group presents the solutions based on the Tsai-Wu failure criterion as a constraint (Eq. 3), where the action of different stress components is considered. It allows obtaining the reliable values of main parameters, applying suggested multi-objective optimization for non-symmetrical composite vessels. These calculations are provided also for three types of composite materials for the studied model. The main goal of the optimization solution is to minimize the amount of composite material and maximize the strength capacity of the vessel. The summary of the main results for the second group is presented in Table. 4. Changes in the main characteristics (thickness, mass, opening radius, and density of strain energy) during the optimization iterative process are presented in Fig. 6-9 for the two groups of solutions.

*Table 4. The summary of the calculations based on the Tsai-Wu criterion*

| Material of composite | Results of calculations | | | |
|---|---|---|---|---|
| | Mass of composite (kg) | | The thickness of the composite on a cylinder (mm) | |
| | initial | final | initial | final |
| **GF-PP** | 26.7 | 20.8 | 10.6 | 8.1 |
| **CF-PA** | 38.8 | 4.5 | 16.8 | 1.9 |
| **FF-PLA** | 27.3 | 16.8 | 11.3 | 7.5 |

The minimization of the total thickness of the composite is illustrated in Fig. 6 for three types of composite materials used for simulations (GF-PP, CF-PA, FF-PLA) for two groups of solutions. This value includes the thicknesses of all helical and hoop layers on the cylinder. Fig. 7a shows the total output composite thickness along the vessel height with a characteristic composite thickness profile with a constant thickness in the cylindrical part and variable thickness in both domes.

One of the key aspects of MIGA methods involves the necessity to set the correct range for parameters. This ensures that, in the process of searching for an optimal solution, the generated values remain within the allowable range of the geometric parameters of the mechanical model. The initial values of all input parameters, including thicknesses, are randomly generated. For each designed layer, the thickness ranges are set from 0.01 to 2 mm, and orientation angles from 0° to 90°. MIGA is applied, with 10 generations, 10 islands, and 10 subpopulations, and solution time for 1000 iterations takes an average of 510 min. This number of iterations (1000 iterations) can be considered sufficient at this stage in our case since the algorithm converged without further improving the minimum feasible values of the objective functions of thickness (Fig. 6) and composite mass (Fig. 7b) at maximum values of strain energy. The comparison of thicknesses, obtained using two types of criteria for each material shows that the difference is lower for the material with better mechanical characteristics (Table 1), such as CF-PA. Test studies were carried out for three types of composites in order to assess the cost-effectiveness of selecting a composite material for the considered water tank. As can be seen from Fig. 6 and Table 4, for CF-PA we obtained thickness and weight 4-4.5 times less compared to other composites (GF-PP, FF-PLA), which is reasonable due to the high carbon material properties.

For each mechanical solution, we obtain the density of strain energy as an output value. The maximum value changes are shown in Fig. 8. The solutions of the first group are characterized by a strong increase in the density of strain energy compared to the solution of the second group. This difference is mainly explained by the fact that initial geometrical parameters and thicknesses for the first group were far from the optimal solution.

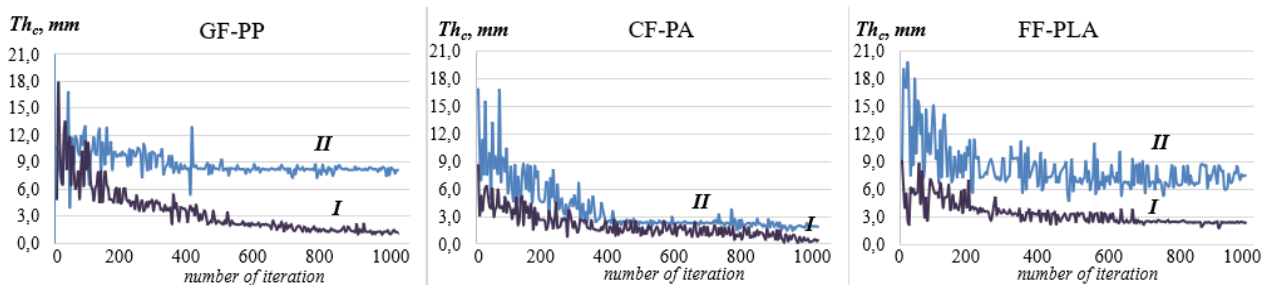

*Fig.6.* *Iterations in total thickness of composite on the cylindrical part during optimization with different failure index: **I** – principal stresses criterion; **II** – Tsai-Wu criterion.*

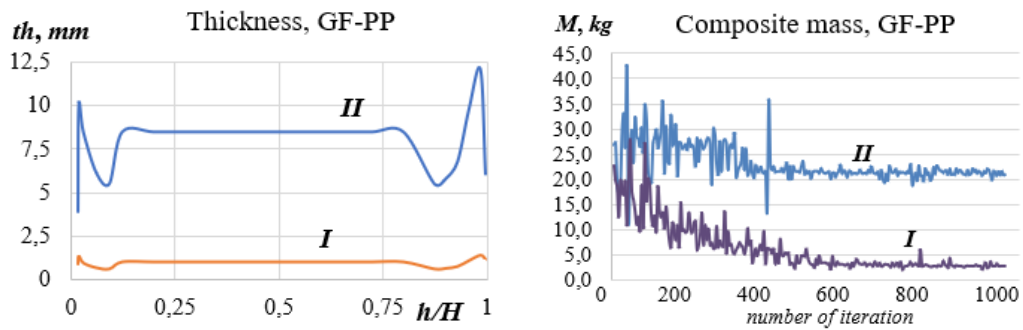

a) b)

*Fig.7. Results of optimization solution for GF-PP with different failure index: **I** – principal stresses criterion; **II** – Tsai-Wu criterion: a) distribution of total composite thickness along the vessel; b) changes of composite mass during optimization.*

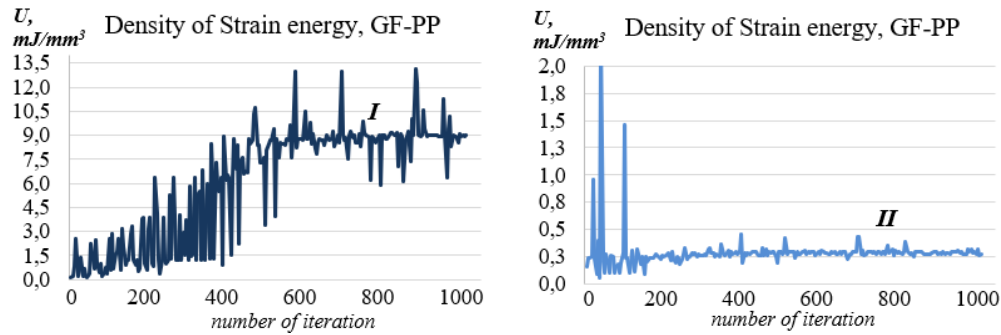

*Fig.8. Changes of the maximal density of strain energy during optimization with different failure index: **I** – principal stresses criterion; **II** – Tsai-Wu criterion for GF-PP.*

The calculations thus point to optimal values for the polar opening radii and, as a consequence, changes in fiber angle for all designed layers in the composite. The changes in polar openings, considering the vessel geometry, are specified for each radius separately. The results are presented for the values of polar opening for upper and lower domes for GF-PP (Fig. 9). Variation during the optimization process lies in the range of feasible values. It can be noted that the radius of the second polar opening on the upper dome $r_{0up_2}$ tends to be the first $r_{0up_1} = 20mm$ (constant) in the iterative process as the thickness of the composite decreases, which provides better strength on the upper dome. The third opening radii $r_{0up_3}$ and $r_{0d_3}$ tend to have approximately the same values on the upper and down domes.

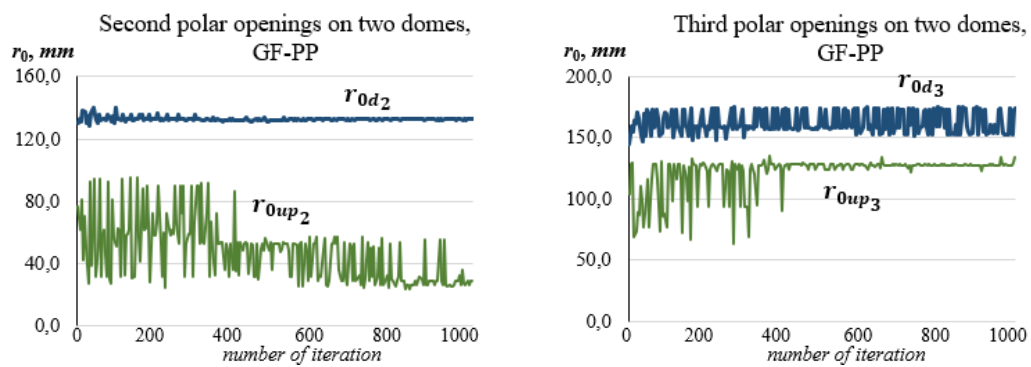

*Fig.9. Changes of second and third polar openings for GF-PP during optimization with Tsai-Wu criterion.*

The developed algorithm allows us to obtain the distribution of fiber angles for each helical paired layer along the height of the vessel, according to the suggested approach. Obtained initial and final distributions are presented in Fig. 10*a,b,c* for GF-PP, CF-PA, FF-PLA, respectively. The distributions of fiber angles for each package of helical layers ($\alpha_1, \alpha_2, \alpha_3$) are fully dependent on the polar openings on each dome ($r_{0up_1}, r_{0up_2}, r_{0up_3}, r_{0d_1}, r_{0d_2}, r_{0d_3}$) (Table 5). $r_{0up_1}$ and $r_{0d_1}$ remain constant during optimization, according to the manufacturing constraints. The zones of the down and upper domes are characterized by a sharp increase in the fiber angle up to 90º at the polar openings. Initial values are estimated, according to the initial random values of polar openings on each dome. It can be noted also from the final redistribution in fiber angles, for each composite material, that the distribution of angles in 3rd layer tends to the constant values on the cylindrical part with equal polar openings on each dome.

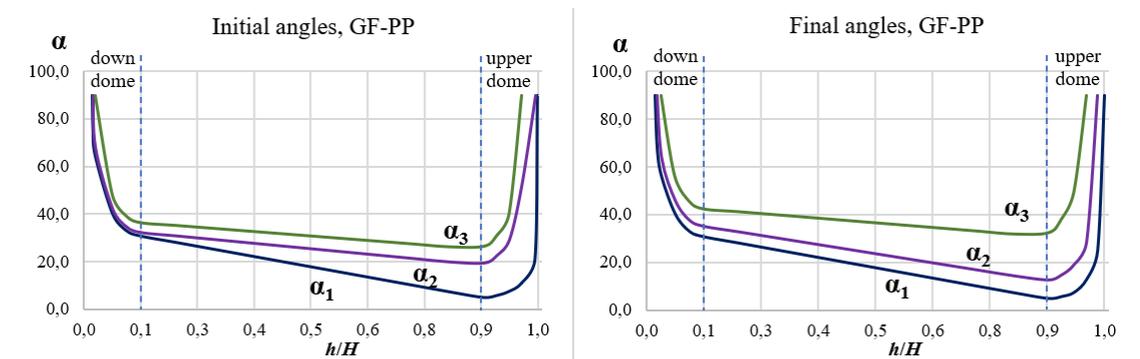

*a)*

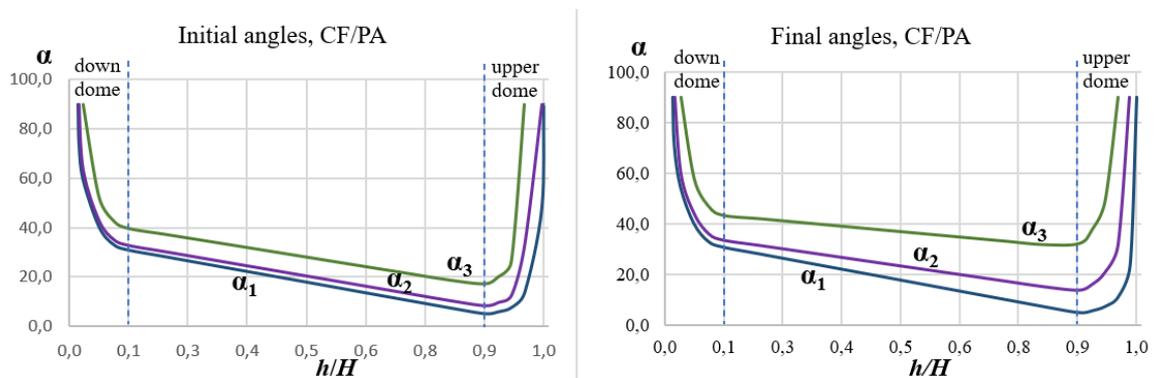

*b)*

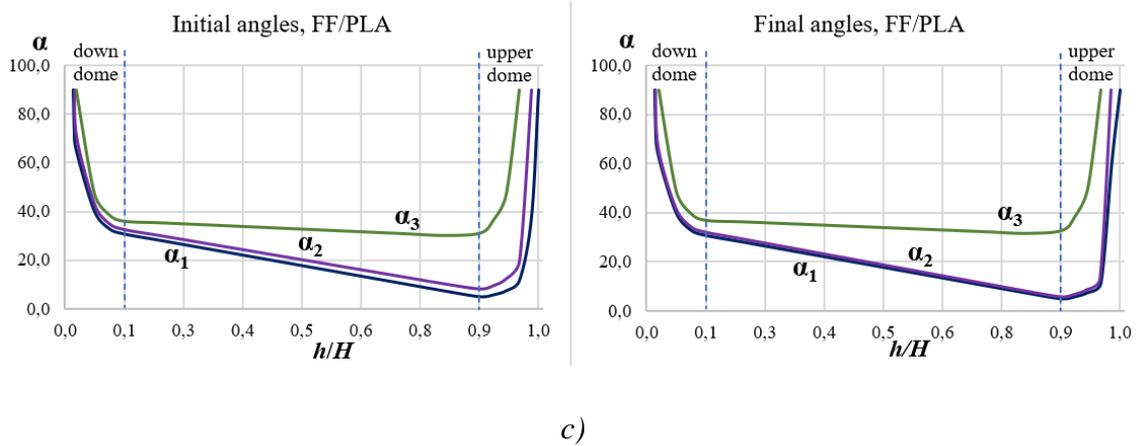

*c)*

*Fig.10. Initial and final distribution of 1st, 2nd, and 3rd helical layers (the calculations based on the Tsai-Wu criterion): a) for GF-PP; b) for CF-PA; c) for FF-PLA.*

*Table 5. Values of polar openings on the upper and down dome for three packages of helical layers, for the second group (the calculations based on the Tsai-Wu criterion)*

| Material/ $r_0$, mm | GF-PP | | CF-PA | | FF-PLA | |
|---|---|---|---|---|---|---|
| | *Initial* | *Final* | *Initial* | *Final* | *Initial* | *Final* |
| $r_{0_1}$ for $\alpha_1$ | | | | | | |
| $r_{0up_1}$=const | 20.0 | --- | 20.0 | --- | 20.0 | --- |
| $r_{0d_1}$= const | 125.0 | --- | 125.0 | --- | 125.0 | --- |
| $r_{0_2}$ for $\alpha_2$ | | | | | | |
| $r_{0up_2}$ | 76.9 | 29.0 | 76.3 | 55.4 | 86.3 | 23.1 |
| $r_{0d_2}$ | 129.9 | 132.0 | 130.3 | 134.8 | 130.4 | 129.7 |
| $r_{0_3}$ for $\alpha_3$ | | | | | | |
| $r_{0up_3}$ | 103.5 | 134.3 | 131.6 | 123.6 | 117.7 | 126.2 |
| $r_{0d_3}$ | 144.8 | 173.6 | 161.4 | 167.3 | 173.0 | 146.0 |

Numerical solutions using the developed methodology also revealed a significant dependency between components of material properties and obtained results. The values of strength components of composite mechanical properties have a notable impact on the resulting optimal thickness, even with just one component. However, accurately assessing their values is often challenging. While Chamis' formulas (Appendix 1, Eq A1) can be utilized, it has been observed that they may provide an overestimated result for the strength components; for instance, considering the value $Y_t$ for GF-PP. The experimental value is 13.1 MPa, while the calculated value is 24 MPa. This component affects the value of the Tsai-Wu criterion (Eq. 3) and, consequently, the final values of composite thickness and mass. The results of calculations in Table 4 are presented using available data for this material, but, the using of the calculated value 24 MPa leads to the decrease of the final optimal thickness value, for example for GF-PP, from 8.1 mm to 6.7 mm (17%). A similar effect can be observed for the values of other components. This highlights the importance of considering this factor when designing a composite vessel using the proposed algorithm. Therefore, the optimization task is highly sensitive to the choice of the failure criterion, and the values of mechanical properties for the composite material are crucial for finding an optimal solution for the composite pressure vessel design.

## 6. Conclusions

In this work, a newly developed iterative algorithm for the design optimization of a composite pressure vessel based on multi-objective parametric optimization is presented. The criterion for choosing optimal parameters is the assessment of the mechanical state of the vessel, fulfilling the failure criterion. The main originality of the developed optimization methodology is that it provides an efficient tool for the design of composite vessels with unequal polar openings. Another originality is the integration of different opening radii for each helical layer as independent optimization parameters, that provide the non-symmetry of a model.

A special Python script was developed that considers manufacturing features of filament winding, such as the distribution of fiber angles for geodesic and non-geodesic paths on each dome and changes in thickness. Also, due to the presence of unequal openings for winding on each dome, an approach was implemented that includes piecewise linear modeling of fiber angle changes along the cylindrical part of the vessel.

The solution to the optimization problem was provided in the optimization module using a genetic algorithm. According to the literature review and own test calculations genetic algorithms allow us to avoid obtaining local solutions for such types of optimization problems, then, for example, gradient or direct search methods. The optimal parameters of a composite vessel were found by minimizing the composite mass and thickness and maximizing the strain energy.

Testing of the developed optimization methodology was carried out for a composite vessel for storage and heating water. Calculations were provided for three types of composite materials (CF, GF, FF) to evaluate the mechanical capability and economic profitability of their use for such a design model of the vessel. Also, the influence of material properties of the composite on the final result was estimated. It was obtained that vessels made of CF-PA can have much lower thickness and weight compared to other composites under consideration (GF-PP, FF-PLA). However, the cost of CF-PA would be much higher. Of course, other aspects must also be considered when choosing a material. The economic study that would follow this work will have to consider the material cost, the mass of materials, but also other criteria such as the recyclability and value of the recyclable material. The considered iterative algorithm determines the optimal vessel parameters, based on the mechanical problem and some design features. For a complete analysis, the heat transfer and thermomechanical problem should also be studied, which may be a topic for further research.

**Acknowledgments**

This work is supported by the Corenstock Industrial Chair ANR-20-CHIN-0004-01, co-founded by Elm.leblanc, IMT and its schools IMT Nord Europe and Mines Saint-Étienne. Funded by the French National Research Agency, the project is also endorsed by the Institut Carnot Télécom & Société Numérique and accredited by the EMC2 cluster.

L. Rozova acknowledges the financial support of a special program of the National Research Agency: PAUSE – ANR Ukraine program: Ukrainian scientists support (ANR-23-PAUK-9035-01) (https://anr.fr/en/call-for-proposals-details/call/pause-anr-ukraine-program-ukrainian-scientists-support/).

# Appendix 1

All elastic and strength properties for the three used materials should be estimated. In this work, the Chamis' micromechanical model was used [38,50].

$$E_{11} = V^f E_{11}^f + V^m E^m$$

$$E_{22} = \frac{E^m}{1 - \sqrt{V^f}\left(1 - E^m/E_{22}^f\right)}$$

$$\nu_{12} = \nu_{13} = V^f \nu_{12}^f + V^m \nu^m$$

$$\nu_{23} = \frac{E_{22}}{2G_{23}} - 1$$

$$G_{12} = \frac{G^m}{1 - \sqrt{V^f}\left(1 - G^m/G_{12}^f\right)}$$

$$G_{23} = \frac{G^m}{1 - \sqrt{V^f}\left(1 - G^m/G_{23}^f\right)} \tag{A1}$$

$$X_t = X_t^f \cdot V^f$$

$$Y_t = \left[1 - \left(\sqrt{V^f} - V^f\right)\left(1 - \frac{E^m}{E_{22}^f}\right)\right] X_t^m$$

$$X_c = X_c^f \cdot V^f$$

$$Y_c = \left[1 - \left(\sqrt{V^f} - V^f\right)\left(1 - \frac{E^m}{E_{22}^f}\right)\right] Y_c^m$$

$$S_{12} = \left[1 - \left(\sqrt{V^f} - V^f\right)\left(1 - \frac{G^m}{G_{12}^f}\right)\right] S_{12}^m$$

where $E_{11}, E_{22}$ - Young's modulus along fiber direction and in transverse direction; $\nu_{12}, \nu_{13}, \nu_{23}$ - Poisson's coefficients, $G_{12}, G_{23}$ - shear modulus. $V$ - component volume fracture, subscripts *f* and *m* stand for fiber and matrix material respectively.

# References

[1] Ahmadi Jebeli M, Heidari-Rarani M. Development of Abaqus WCM plugin for progressive failure analysis of type IV composite pressure vessels based on Puck failure criterion. Eng Fail Anal 2022;131:105851. https://doi.org/10.1016/j.engfailanal.2021.105851.

[2] Regassa Y, Gari J, Lemu HG. Composite Overwrapped Pressure Vessel Design Optimization Using Numerical Method. Journal of Composites Science 2022;6:229. https://doi.org/10.3390/jcs6080229.

[3] Zhou W, Wang J, Pan Z, Liu J, Ma L, Zhou J, et al. Review on optimization design, failure analysis and non-destructive testing of composite hydrogen storage vessel. Int J Hydrogen Energy 2022;47:38862–83. https://doi.org/10.1016/j.ijhydene.2022.09.028.

[4] Forth SC., Pat B. Composite Overwrapped Pressure Vessels, a Primer; NASA Center for Aero Space Information: Hanover, MD, USA, 2011 n.d.

[5] Solazzi L, Vaccari M. Reliability design of a pressure vessel made of composite materials. Compos Struct 2022;279:114726. https://doi.org/10.1016/j.compstruct.2021.114726.

[6] Alcántar V, Aceves SM, Ledesma E, Ledesma S, Aguilera E. Optimization of Type 4 composite pressure vessels using genetic algorithms and simulated annealing. Int J Hydrogen Energy 2017;42:15770–81. https://doi.org/10.1016/j.ijhydene.2017.03.032.

[7] Christensen P.W., Klarbring A. An Introduction to Structural Optimization. SpriSc&BusMedia: Berlin/Heildeberg, Germany 2008;153.

[8] Boyd S, Vandenberghe L. Convex Optimization. Cambridge University Press; 2004. https://doi.org/10.1017/CBO9780511804441.

[9] Janga Reddy M, Nagesh Kumar D. Evolutionary algorithms, swarm intelligence methods, and their applications in water resources engineering: a state-of-the-art review. H2Open Journal 2020;3:135–88. https://doi.org/10.2166/h2oj.2020.128.

[10] Mei L, Wang Q. Structural Optimization in Civil Engineering: A Literature Review. Buildings 2021;11:66. https://doi.org/10.3390/buildings11020066.

[11] Dassault Systèmes D. Abaqus analysis user's guide. Technical Report Abaqus 6.14 Documentation, Simulia Corp 2016.

[12] Rafiee R, Shahzadi R, Speresp H. Structural optimization of filament wound composite pipes. Frontiers of Structural and Civil Engineering 2022;16:1056–69. https://doi.org/10.1007/s11709-022-0868-3.

[13] Tyflopoulos E, Hofset TA, Olsen A, Steinert M. Simulation-based design: a case study in combining optimization methodologies for angle-ply composite laminates. Procedia CIRP 2021;100:607–12. https://doi.org/10.1016/j.procir.2021.05.131.

[14] Mian HH, Wang G, Dar UA, Zhang W. Optimization of Composite Material System and Lay-up to Achieve Minimum Weight Pressure Vessel. Applied Composite Materials 2013;20:873–89. https://doi.org/10.1007/s10443-012-9305-4.

[15] Veivers H, Bermingham M, Dunn M, Veidt M. Layup optimisation of laminated composite tubular structures under thermomechanical loading conditions using PSO. Compos Struct 2021;276:114483. https://doi.org/10.1016/j.compstruct.2021.114483.

[16] Almeida JHS, St-Pierre L, Wang Z, Ribeiro ML, Tita V, Amico SC, et al. Design, modeling, optimization, manufacturing and testing of variable-angle filament-wound cylinders. Compos B Eng 2021;225:109224. https://doi.org/10.1016/j.compositesb.2021.109224.

[17] Daghighi S, Zucco G, Weaver PM. Design of Variable Stiffness Super Ellipsoidal Pressure Vessels under Thermo-mechanical Loading. AIAA SCITECH 2022 Forum, Reston, Virginia: American Institute of Aeronautics and Astronautics; 2022. https://doi.org/10.2514/6.2022-0869.

[18] Rafiee R, Shahzadi R, Jafari S. Filament wound pipes optimization platform development: A methodological approach. Compos Struct 2022;297:115972. https://doi.org/10.1016/j.compstruct.2022.115972.

[19] Daghighi S, Weaver PM. Nonconventional tow-steered pressure vessels for hydrogen storage. Compos Struct 2024;334:117970. https://doi.org/10.1016/j.compstruct.2024.117970.

[20] Kang C, Liu Z, Shirinzadeh B, Zhou H, Shi Y, Yu T, et al. Parametric optimization for multi-layered filament-wound cylinder based on hybrid method of GA-PSO coupled with local sensitivity analysis. Compos Struct 2021;267:113861. https://doi.org/10.1016/j.compstruct.2021.113861.

[21] Kim C-U, Hong C-S, Kim C-G, Kim J-Y. Optimal design of filament wound type 3 tanks under internal pressure using a modified genetic algorithm. Compos Struct 2005;71:16–25. https://doi.org/10.1016/j.compstruct.2004.09.006.

[22] Ellul B, Camilleri D. The applicability and implementation of the discrete Big Bang-Big Crunch optimisation technique for discontinuous objective function in multi-material laminated composite pressure vessels. International Journal of Pressure Vessels and Piping 2018;168:39–48. https://doi.org/10.1016/j.ijpvp.2018.08.008.

[23] Cvetkoska D, Dimovski I, Samak S, Trompeska M, Dukovski V. Using Constrained Multi-Optimization in Design of Composite for Filament Wound High Pressure Vessels. International Journal of Mathematics Trends and Technology 2018;61:107–16. https://doi.org/10.14445/22315373/IJMTT-V61P516.

[24] Qiu C, Han Y, Shanmugam L, Zhao Y, Dong S, Du S, et al. A deep learning-based composite design strategy for efficient selection of material and layup sequences from a given database. Compos Sci Technol 2022;230:109154. https://doi.org/10.1016/j.compscitech.2021.109154.

[25] Ehsani A, Dalir H. Multi-objective optimization of composite angle grid plates for maximum buckling load and minimum weight using genetic algorithms and neural networks. Compos Struct 2019;229:111450. https://doi.org/10.1016/j.compstruct.2019.111450.

[26] Liu X, Qin J, Zhao K, Featherston CA, Kennedy D, Jing Y, et al. Design optimization of laminated composite structures using artificial neural network and genetic algorithm. Compos Struct 2023;305:116500. https://doi.org/10.1016/j.compstruct.2022.116500.

[27] Zu L, Koussios S, Beukers A. Design of filament–wound domes based on continuum theory and non-geodesic roving trajectories. Compos Part A Appl Sci Manuf 2010;41:1312–20. https://doi.org/10.1016/j.compositesa.2010.05.015.

[28] Zu L, Koussios S, Beukers A. Shape optimization of filament wound articulated pressure vessels based on non-geodesic trajectories. Compos Struct 2010;92:339–46. https://doi.org/10.1016/j.compstruct.2009.08.013.

[29] Liang C-C, Chen H-W, Wang C-H. Optimum design of dome contour for filament-wound composite pressure vessels based on a shape factor. Compos Struct 2002;58:469–82. https://doi.org/10.1016/S0263-8223(02)00136-8.

[30] Daghighi S, Rouhi M, Zucco G, Weaver PM. Bend-free design of ellipsoids of revolution using variable stiffness composites. Compos Struct 2020;233:111630. https://doi.org/10.1016/j.compstruct.2019.111630.

[31] Guo K, Wen L, Xiao J, Lei M, Wang S, Zhang C, et al. Design of winding pattern of filament-wound composite pressure vessel with unequal openings based on non-geodesics. J Eng Fiber Fabr 2020;15:155892502033397. https://doi.org/10.1177/1558925020933976.

[32] Li H, Liang Y, Bao H. Splines in the parameter domain of surfaces and their application in filament winding. Computer-Aided Design 2007;39:268–75. https://doi.org/10.1016/j.cad.2006.12.003.

[33] Zu L, Koussios S, Beukers A. Design of filament-wound isotensoid pressure vessels with unequal polar openings. Compos Struct 2010;92:2307–13. https://doi.org/10.1016/j.compstruct.2009.07.013.

[34] Vargas Rojas E, Chapelle D, Perreux D, Delobelle B, Thiebaud F. Unified approach of filament winding applied to complex shape mandrels. Compos Struct 2014;116:805–13. https://doi.org/10.1016/j.compstruct.2014.06.009.

[35] Azeem M, Ya HH, Alam MA, Kumar M, Stabla P, Smolnicki M, et al. Application of Filament Winding Technology in Composite Pressure Vessels and Challenges: A Review. J Energy Storage 2022;49:103468. https://doi.org/10.1016/j.est.2021.103468.

[36] Park CH, Lee W Il. Manufacturing: Economic Consideration. Wiley Encyclopedia of Composites, Wiley; 2012, p. 1–13. https://doi.org/10.1002/9781118097298.weoc129.

[37] Vasiliev V V., Morozov E V. Composite Pressure Vessels. Advanced Mechanics of Composite Materials and Structures, Elsevier; 2018, p. 787–821. https://doi.org/10.1016/B978-0-08-102209-2.00012-8.

[38] Hu N, editor. Composites and Their Properties. InTech; 2012. https://doi.org/10.5772/2816.

[39] De Luca A, Caputo F. A review on analytical failure criteria for composite materials. AIMS Mater Sci 2017;4:1165–85. https://doi.org/10.3934/matersci.2017.5.1165.

[40] Rafiee R, Fakoor M, Hesamsadat H. The influence of production inconsistencies on the functional failure of GRP pipes. Steel and Composite Structures 2015;19:1369–79. https://doi.org/10.12989/scs.2015.19.6.1369.

[41] Daghighi S, Weaver PM. Three-dimensional effects influencing failure in bend-free, variable stiffness composite pressure vessels. Compos Struct 2021;262:113346. https://doi.org/10.1016/j.compstruct.2020.113346.

[42] Tew BW. Preliminary Design of Tubular Composite Structures Using Netting Theory and Composite Degradation Factors. J Press Vessel Technol 1995;117:390–4. https://doi.org/10.1115/1.2842141.

[43] https://www.avient.com/sites/default/files/2020-08/polystrand-selection-guide-1.pdf.

[44] https://www.toraytac.com/product-explorer/products/r2Vf/Toray-Cetex-TC910

[45] https://www.groupedepestele.com/pro_ecomateriaux_lincoreplaff.html

[46] Zhu B, Yu TX, Tao XM. Large Shear Deformation of E-glass/ Polypropylene Woven Fabric Composites at Elevated Temperatures. Journal of Reinforced Plastics and Composites 2009;28:2615–30. https://doi.org/10.1177/0731684408093095.

[47] Çuvalci H, Erbay K, İpek H. Investigation of the Effect of Glass Fiber Content on the Mechanical Properties of Cast Polyamide. Arab J Sci Eng 2014;39:9049–56. https://doi.org/10.1007/s13369-014-1409-8.

[48] Weber DE, Graupner N, Müssig J. Manufacturing of flax- and glass-fibre reinforced thin-walled tubes and measuring their interlaminar shear properties by torsion tests. Compos Struct 2023;319:117191. https://doi.org/10.1016/j.compstruct.2023.117191.

[49] https://www.3ds.com/products-services/simulia/products/isight-simulia-execution-engine

[50] Chamis C. Mechanics of Composite Materials: Past, Present and Future. Journal of Composites Technology & Research 1989:3–14.